\documentclass[]{imsart}

\RequirePackage[colorlinks,citecolor=blue,urlcolor=blue]{hyperref}

\usepackage{graphicx,color,enumerate,amssymb}
\usepackage{amsmath}
\usepackage{epsfig}
\usepackage{mathrsfs}
\usepackage{multirow}
\usepackage{lscape}

\startlocaldefs

\def\N{{\rm I\kern-.15em N}}
\def\R{{\rm I\kern-.2em R}}
\def\Z{{\rm Z\kern-.26em Z}}

\newtheorem{thm}{Theorem}[section]

\newtheorem{rem}[thm]{Remark}
\newtheorem{cor}[thm]{Corollary}

\newcommand{\be}{\begin{eqnarray}}
\newcommand{\ee}{\end{eqnarray}}
\newcommand{\bq}{\begin{eqnarray*}}
\newcommand{\eq}{\end{eqnarray*}}

\newcommand{\verk}{\stackrel{{\cal D}}{\longrightarrow}}

\newcommand{\bewend}{\hspace*{2mm}\rule{3mm}{3mm}}

\newcommand{\RR}{\mathbb{R}}
\newcommand{\BE}{\mathbb{E}}
\newcommand{\PP}{\mathbb{P}}
\newcommand{\HH}{\mathbb{H}}
\def\L{\text{L}}
\def\C{\text{C}}
\def\S{\text{S}}

\endlocaldefs

\begin{document}

\begin{frontmatter}

\title{A unified approach to goodness-of-fit testing for spherical and hyperspherical data}
\runtitle{Goodness-of-fit testing for hyperspherical data}

 \author{\fnms{Bruno} \snm{Ebner}\corref{}\ead[label=e1]{bruno.ebner@kit.edu}\ead[label=e4,url]{www.math.kit.edu/stoch/~ebner/}}
 \address{Institute of Stochastics, Karlsruhe Institute of Technology (KIT) \\ Englerstr. 2, 76131 Karlsruhe, Germany \\ \printead{e1}\\ \printead{e4}}
\and
\author{\fnms{Norbert} \snm{Henze}\ead[label=e2]{norbert.henze@kit.edu}}
\address{Institute of Stochastics, Karlsruhe Institute of Technology (KIT) \\ Englerstr. 2, 76131 Karlsruhe, Germany \\ \printead{e2}}
\and \author{\fnms{Simos} \snm{Meintanis}\ead[label=e3]{simosmei@econ.uoa.gr}}
\address{Department of Economics, National and Kapodistrian University
of Athens, Athens, Greece}\address{Pure and Applied Analytics, North--West University, Potchefstroom, South Africa\footnote{On sabbatical leave from the University of Athens}\\ \printead{e3}}

\runauthor{B. Ebner, N. Henze and S. Meintanis}

\begin{abstract}
We propose a general and relatively simple method for the construction of goodness-of-fit tests on the sphere and the hypersphere. The method is based on the characterization of probability distributions via their characteristic function, and it leads to test criteria that are convenient regarding applications  and consistent against arbitrary deviations from the model under test. We emphasize goodness-of-fit tests for spherical distributions due to their importance in applications and  the  relative scarcity of available methods.
\end{abstract}





\end{frontmatter}

\section{Introduction} \label{intro}
Let $\|\cdot \|$ be the Euclidean norm in $\RR^d$, $d \ge 2$, and write ${\cal S}^{d-1} := \{x \in \RR^d: \|x\|=1\}$ for the
surface of the unit sphere in $\RR^d$. In this paper, we consider the problem of testing goodness-of-fit for distributions defined on the sphere ${\cal S}^{2}$ or on the hypersphere ${\cal S}^{d-1}$, where $d>3$. In this respect, there is a plethora of such tests for distributions defined on $\mathbb R^d$, even in the multivariate case $d>1$. Besides, also goodness-of-fit tests on the circular domain ${\cal S}^{1}$
is a relatively well-explored area. For the latter case we refer, e.g., to \cite{GV:2018} (Chapters 6 and 7), \cite{JJM:2019} and \cite{JMV:2020}.

On the other hand, the same problem for data taking values on ${\cal S}^{d-1}$, where $d\geq 3$, has been mostly confined to testing for uniformity. Nevertheless, and while  the notion of ``non-preferred direction" and hence testing for uniformity is certainly central to (hyper)spherical data analysis, there are several more flexible distributions, which in fact often have the uniform as a special case. The reader is referred to the monographs of \cite{LV:2017}, Section 2.3, and \cite{MJ:2000}, Section 9.3, for such non-uniform models for (hyper)spherical data. At the same time, it seems  that  goodness-of-fit tests specifically tailored to hyper(spherical) laws are scarce, certainly in the case of a composite null hypothesis, where distributional parameters need to be estimated from the data at hand, but also for a completely specified hypothesis with fixed (known) parameter values. For the latter case, the test based on nearest neighbors  proposed in \cite{EHY:2018} seems to be one of the few tests available, while to the best of our knowledge, there is much need for research in the case of a composite hypothesis.

In view of these lines, we suggest a procedure for testing goodness-of-fit for distributions defined on ${\cal S}^{d-1}$, where $d\geq 3$.\footnote{The new test also applies to circular distributions ($d=2$), but herein we emphasize the higher dimensional cases.} The suggested  test is novel in that it is general-purpose  suitable for arbitrary (hyper)spherical distributions, either with fixed or estimated parameters,  and it is straightforwardly applicable  provided that one can easily draw Monte Carlo samples from the distribution under test.

Suppose $X$ is a random (column) vector in $\RR^d$ taking  values in ${\cal S}^{d-1}$
with a density $f$ with respect to surface measure and characteristic function (CF) $\varphi(t)=\BE ({\rm e}^{{\rm{i}}t^\top X})$,  $t \in \RR^{d}$, where $\top$ denotes transpose,
and ${\rm{i}}=\sqrt{-1}$ stands for the imaginary unit.
 \smallskip
We start our exposition with the simple null hypothesis
\be \label{null}
H_0: f=f_0,
\ee
where $f_0$ is some given density on ${\cal S}^{d-1}$, which should be tested against
the general alternative $H_A$ that the distributions pertaining to $f$ and $f_0$ are different.  If $X_0$ has density $f_0$ and CF
$\varphi_0(t) = \BE ({\rm e}^{{\rm{i}}t^\top X_0})$,  $t \in \RR^{d}$, say, the standard CF-based statistic for testing $H_0$ versus $H_A$
is given by
\begin{equation} \label{tw1}
D_{n,w}= \int |\varphi_n(t)-\varphi_0(t)|^2 w(t) \, {\rm{d}} t.
\end{equation}
Here,
\be \label{CF}
\varphi_n(t)=\frac{1}{n}\sum_{j=1}^n\exp({\rm i} t^\top  X_j), \ t \in \RR^d,
\ee
is the empirical CF of $X_1,\ldots,X_n$, and $X_1,\ldots,X_n$ are independent and identically distributed (i.i.d.) copies of $X$.
In \eqref{tw1}, the domain of integration as well as the nonnegative weight function $w(\cdot)$ will be specified below
in a way that $D_{n,w}$ is amenable  to computation and that a test of $H_0$ that rejects $H_0$ for large values of $D_{n,w}$ is consistent
against each alternative to $H_0$. Notice that $D_{n,w}$ is an estimator of the
population Fourier-type discrepancy measure
\begin{equation} \label{tw}
D_{w}(f,f_0)= \int |\varphi(t)-\varphi_0(t)|^2 w(t) \, {\rm{d}} t
\end{equation}
between $f$ and $f_0$.
\medskip

The starting point of this paper is that the approach outlined above assumes that the  functional form
of the CF $\varphi_0$ is known. Such knowledge, however, is only available for distributions on the real line $\mathbb R^1$ and for a few selected cases of multivariate distributions, such as the multivariate normal and the multivariate stable distribution; see \cite{EH:2020} and \cite{MNT:2015}.

In order to circumvent this obstacle, which is even more challenging for distributions taking values on ${\cal{S}}^{d-1}$, we suggest the test statistic
\be \label{ts}
T_{n,m,w} = \frac{mn}{m+n} \int |\varphi_n(t)-\psi_{m}(t)|^2 w(t)\, {\rm{d}} t.
\ee
Here,
\be \label{CFnull}
\psi_{m}(t)=\frac{1}{m}\sum_{j=1}^m\exp({\rm i} t^\top Y_{j}), \quad t \in \RR^d,
\ee
is the empirical CF of $Y_{1}, \ldots, Y_{m}$, where, independently of $X_1,\ldots,X_n$, the random vectors $Y_{1}, \ldots, Y_{m}$ are i.i.d. copies of $X_0$.
Of course, realizations of $Y_{1}, \ldots, Y_{m}$ are generated via Monte Carlo.
Notice that $\psi_{m}$ is an estimator of the CF $\varphi_{0}$. In this way, the functional form of $\varphi_0$
is not needed in the test statistic $T_{n,m,w}$ in \eqref{ts}, which is reminiscent of a CF-based test for the two-sample problem,
 one sample being the data $X_1,\ldots,X_n$ at hand, while the other consists of artificial data
 generated under the null hypothesis $H_0$.
  For more details on CF-based tests for the two-sample problem the reader is referred to \cite{M:2005} and \cite{AJG:2008}.
\medskip

This idea also applies to the problem of testing  the composite  null hypothesis
 \be \label{null2} H_{0,\vartheta}: f(\cdot) =f_0(\cdot ,\vartheta) \quad \text{for some} \ \vartheta \in \Theta,
\ee
 against general alternatives.  Here, $\{f_0(\cdot,\vartheta): \vartheta \in \Theta\}$
 is a given family of densities on ${\cal{S}}^{d-1}$  that is parameterized  in terms of $\vartheta \in \Theta$, where $\Theta \subset \RR^s$ for some $s \ge 1$.
 In this setting, the test statistic in \eqref{ts} is modified according to
 \be \label{ts3}
\widehat T_{n,m,w}= \frac{mn}{m+n} \int |\varphi_n(t)-\widehat \psi_{m}(t)|^2 w(t)\, {\rm{d}}t,
\ee
with $\varphi_n(t)$  defined in \eqref{CF} and
\be \label{CF3}
\widehat \psi_{m}(t)=\frac{1}{m}\sum_{j=1}^m\exp({\rm i}t^\top \widehat {Y}_{j}),\quad t \in \RR^d,
\ee
where $\widehat {Y}_{j}$, $j \in \{1,...,m\}$, are i.i.d. copies of a random vector having density $f_0(\cdot;\widehat {\vartheta}_n)$.
Here, $ \widehat {\vartheta}_n:=\widehat {\vartheta}_n(X_1,...,X_n)$ is some estimator of $\vartheta$ computed from $X_1,\ldots,X_n$.
\bigskip

 In this connection, we note that the idea of a goodness-of-fit method that employs an artificial sample from the distribution under test seems to date back to \cite{friedman03}, at least for independent data and simple hypotheses. Recently, \cite{chen22} proposed a CF-based method using the notion of artificial samples for goodness-of-fit within the family of multivariate elliptical distributions, \cite{CX:2023}
 employ artificial samples in order to specifically test multivariate normality in high dimensions using nearest neighbors, while
 \cite{salmaso} applies a test procedure  for mixed data by means of artificial samples.

The remainder of this work unfolds as follows. In Section~\ref{limits} we obtain the limit null distribution of $T_{n,m,w}$ as well as the corresponding law under fixed deviations from $H_0$. In Section \ref{boot}, the validity of a bootstrap resampling scheme necessary for actually carrying out the test  for simple hypotheses with fixed parameters is established, while in Section~\ref{boot1} a corresponding bootstrap resampling for the composite hypothesis test statistic $\widehat T_{n,m,w}$ is suggested. Section \ref{simul} contains an extensive Monte Carlo study of the finite-sample behavior of the new tests including comparisons, while Section~\ref{real} illustrates real-data applications. The final Section~\ref{conclude} provides some discussion.

\section{Asymptotics} \label{limits}
In this section, we provide the limit distribution of $T_{n,m,w}$ defined in \eqref{ts}. To be flexible
with respect to both the region of integration and to the weight function $w$, let $M$ be some nonempty
Borel set in $\RR^d$, and let $\mu$ be some finite measure on (the Borel subsets of) $M$. Thus, $M$ could  be $\RR^d$ itself, and $\mu$ could be absolutely continuous with respect to the Lebesgue measure in $\RR^d$,  or $M$ could
be ${\cal S}^{d-1}$, and $\mu$ could be absolutely continuous with respect to spherical measure. Notably, $M$ could also be some countable subset $T$ of $\RR^d$,
with $\mu$ having a probability mass with respect to the counting measure on $T$.\footnote{A counting measure would be sufficient for circular distributions; see for instance \cite{JJM:2019}.}
\smallskip

In this setting, let $X,X_1,X_2, \ldots$ and $Y,Y_1,Y_2, \ldots $ be independent $M$-valued random vectors that are defined on
some common probability space $(\Omega,{\cal A},\PP)$. Moreover, let $X,X_1,X_2, \ldots$ be i.i.d. with density $f$ with respect to $\mu$ and
CF $\varphi(t) = \BE [\exp({\rm i} t^\top X)]$, $t \in \RR^d$. Furthermore, let $Y,Y_1,Y_2, \ldots $ be i.i.d. with density $g$ with respect to $\mu$ and
CF $\psi(t) = \BE [\exp({\rm i} t^\top Y)]$, $t \in \RR^d$. Recall that
\[
\varphi_n(t) := \frac{1}{n}\sum_{j=1}^n\exp({\rm i} t^\top  X_j), \qquad \psi_m(t)=\frac{1}{m}\sum_{j=1}^m\exp({\rm i} t^\top  Y_j), \quad t \in \RR^d,
\]
are the empirical CF's of $X_1,\ldots,X_n$ and $Y_1,\ldots,Y_m$, respectively. This section tackles the limit distribution of
\begin{equation*}\label{seca-tnm}
T_{n,m} := \frac{nm}{n+m} \int_M \big{|} \varphi_n(t) - \psi_m(t)\big{|}^2 \, \mu(\text{d}t)
\end{equation*}
as $m,n \to \infty$, under each of the conditions $\varphi = \psi$ and $\varphi \neq \psi$.
\medskip

Putting
\[
\C(x) := \int_M \cos\big(t^\top x\big) \, \mu(\text{d}t), \quad  \S(x) := \int_M \sin\big(t^\top x\big) \, \mu(\text{d}t), \quad x \in \RR^d,
\]
straightforward algebra yields
\begin{eqnarray*}
\frac{m+n}{mn} T_{n,m} & = & \frac{1}{n^2}\sum_{j,k=1}^n \big( \C(X_k-X_j) + {\rm i}\, \S(X_k-X_j)\big)\\
& & - \frac{1}{mn} \sum_{k=1}^m \sum_{j=1}^n \big( \C(X_j-Y_k) + {\rm i}\, \S(X_j-Y_k)\big)\\
& & - \frac{1}{mn} \sum_{k=1}^m \sum_{j=1}^n \big( \C(Y_k-X_j) + {\rm i}\, \S(Y_k-X_j)\big)\\
& & + \frac{1}{m^2} \sum_{j,k =1}^m  \big(\C(Y_k-Y_j) + {\rm i}\, \S(Y_k-Y_j)\big).
\end{eqnarray*}
Since $\S(-x) = - \S(x)$ and $\S(0_d) =0$, where $0_d$ is the origin in $\RR^d$, the sum of the imaginary parts vanishes, and we obtain
\begin{equation}\label{eq:test}
\frac{m+n}{mn} T_{n,m}  =  \frac{1}{n^2} \sum_{j,k=1}^n  \C(X_j-X_k) - \frac{2}{mn} \sum_{j=1}^n \sum_{k=1}^m  \C(X_j-Y_k)
 + \frac{1}{m^2} \sum_{j,k=1}^m \C(Y_j-Y_k).
\end{equation}
A further simplification is obtained if we assume that the set $M$ --  like $\mathbb R^{d}$,  ${\cal S}^{d-1}$ or the grid $\mathbb{Z}^d$ --  is symmetric
with respect to the origin $0_d$, i.e., we have $-M =M$, where $-M := \{-x: x \in M\}$. Furthermore, we suppose that the measure $\mu$ is invariant
with respect to the reflection $T(x) := -x$, $x \in \RR^d$, i.e., we have $\mu = \mu^T$, where $\mu^T$ is the image of $\mu$ under $T$. By transformation
of integrals, we then obtain
\[
\S(x) = \int_M \sin(t^\top x) \, \mu(\text{d}t) = - \S(x), \quad x \in \RR^d,
\]
and thus $\S(x) = 0$, $x \in \RR^d$. Putting
\[
\text{CS}(\xi) := \cos \xi  + \sin \xi,\quad \xi \in \RR,
\]
and using the addition theorem $\cos(\alpha - \beta) = \cos  \alpha \cos  \beta + \sin \alpha \sin  \beta$,
some algebra yields
\begin{equation}\label{reprtmn}
T_{n,m} = \frac{mn}{m+n} \int_M \left(\frac{1}{n}\sum_{j=1}^n \text{CS}(t^\top X_j) - \frac{1}{m} \sum_{k=1}^m \text{CS}(t^\top Y_k)\right)^2 \mu(\text{d}t).
\end{equation}
Now, writing ${\cal B}(M)$ for the $\sigma$-field of Borel sets on $M$, let $\HH := \L^2(M,{\cal B}(M),\mu)$ be the separable Hilbert space of (equivalence classes of) measurable functions $u:M \rightarrow \RR$
satisfying $\int_M u^2 \, \text{d} \mu  < \infty$, equipped   with the inner product $\langle u,v\rangle = \int_M u v \, \text{d} \mu$ and the norm
$\|u\|_\HH = \langle u,u \rangle^{1/2}$, $u \in \HH$.

\medskip

\begin{thm}\label{thm1} Suppose that $\varphi = \psi$. If $M$ is symmetric with respect to $0_d$ and $\mu$ is invariant with respect to reflections
at $0_d$, there is a centred Gaussian random element $W$ of $\HH$ having covariance kernel $K(s,t) = \BE\big[ W(s) W(t)\big]$, where
\begin{equation}\label{defcovk}
K(s,t) = {{\rm {Cov}}}\big({\rm{CS}}(s^\top X),{\rm{CS}}(t^\top X)\big), \quad s,t \in M,
\end{equation}
such that $T_{n,m} \verk \|W\|^2_\HH $  as $n,m \to \infty$.
\end{thm}
\smallskip

{\sc Proof.} Let $\varrho(t) := \BE[\text{CS}(t^\top X)]$, $t \in \RR^d$. Notice that $\varrho(t) = \BE[\text{CS}(t^\top Y)]$ since $\varphi = \psi$. Let
\begin{equation}\label{defamnbmn}
a_{m,n} = \sqrt{\frac{m}{m+n}}, \quad b_{m,n} = \sqrt{\frac{n}{m+n}},
\end{equation}
and put
\begin{equation}\label{defunvm}
U_n(t) = \frac{1}{\sqrt{n}} \sum_{j=1}^n \big\{ \text{CS}(t^\top X_j)\! - \!  \varrho(t) \big\}, \quad
V_m(t) = \frac{1}{\sqrt{m}} \sum_{k=1}^m \big\{ \text{CS}(t^\top Y_k)\! - \! \varrho(t) \big\},
\end{equation}
$t \in M$. Then $U_n = U_n(\cdot)$ and $V_m = V_m(\cdot)$ are centred random elements of $\HH$, and by the
central limit theorem in separable Hilbert spaces (see, e.g., Theorem 2.7 of \cite{B:2000}),  we have $U_n \verk U$ as $n \to \infty$ and $V_m \verk V$ as
$m \to \infty$, where $U$ and $V$ are centred random elements of $\HH$ having  covariance kernel $K$  given in \eqref{defcovk}.
From \eqref{reprtmn}, we obtain
\begin{equation}\label{asdisttnm}
T_{n,m} = \big{\|} W_{n,m}\big{\|}^2_\HH ,
\end{equation}
where $W_{n,m} = a_{m,n} U_n - b_{m,n}V_m$. Since $U_n$ and $V_m$ are independent for each pair $(n,m)$, also $U$ and $V$ are independent,
and we have $(U_n,V_m) \verk (U,V)$ as $n,m \to \infty$ (see, e.g., Theorem 2.8. of \cite{B:1999}). Notice that $n,m \to \infty$ means that, as $s \to \infty$, we have
$n= n(s) \to \infty$ and $m = m(s) \to \infty$ in an arbitrary manner. Now, if some subsequence of the bounded sequence $(a_{m,n})$ converges to
some $\tau \in [0,1]$ then, because of the continuous mapping theorem and $a_{m,n}^2 + b_{m,n}^2 =1$, it follows that
$W_{n,m} \verk \sqrt{\tau} U - \sqrt{1-\tau} V$ as $n,m \to \infty$. This limit random element has the same distribution as $W$, irrespective of
$\tau$. Consequently, we have $W_{n,m} \verk W$ as $n,m \to \infty$ (see, e.g., Theorem 2.6 of \cite{B:1999}).  In view of \eqref{asdisttnm} and the continuous mapping theorem, the assertion  follows. \bewend
\bigskip

The next result gives the almost sure limit of $(m+n)T_{n,m}/(mn)$ as $m,n \to \infty$.


\begin{thm}\label{thm2} Let
\begin{equation}\label{defatbt}
\varrho(t) := \BE[{\rm{CS}}(t^\top X)], \quad b(t) := \BE[{\rm{CS}}(t^\top Y)], \quad t \in M.
\end{equation}
 Under the conditions
on $M$ and $\mu$ stated in Theorem~\ref{thm1}, we have
\[
\lim_{n,m \to \infty} \frac{m+n}{mn} T_{n,m}  = \int_M \big(\varrho(t)-b(t)\big)^2 \, \mu({\rm{d}}t) \quad \PP\text{-almost  surely}.
\]
\end{thm}
\smallskip

{\sc Proof.} Let
\begin{equation}\label{defanbm}
A_n(t) = \frac{1}{n} \sum_{j=1}^n \text{CS}(t^\top X_j), \quad B_m(t) = \frac{1}{m} \sum_{k=1}^m \text{CS}(t^\top Y_k)
\end{equation}
and, regarded as random elements of $\HH$,  put $A_n =  A_n(\cdot)$ and $B_m = B_m(\cdot)$. Likewise,
write $a = \varrho(\cdot)$ and $b = b(\cdot)$ for the degenerate random elements of $\HH$ that are the expectations
of $A_n$ and $B_m$, respectively. By the strong law of large numbers  in Banach spaces (see, e.g., \cite{HP:1976}), we have
$\|A_n - a\|_\HH \to 0$ as $n \to \infty$ and $\|B_m - b \|_\HH \to 0$ as $m \to \infty$ $\PP$-almost surely. It follows that
\[
\frac{m+n}{mn} T_{n,m} = \|A_n - B_m\|^2_\HH \to \|a-b\|^2_\HH \quad \PP\text{-almost surely as } n,m \to \infty. \bewend
\]

\medskip
Let
\begin{equation*}\label{deltadef}
\Delta := \int_M |\varphi(t) - \psi(t)|^2\, \mu(\text{d}t).
\end{equation*}
It is readily seen that $\Delta$ equals the almost sure limit figuring in  Theorem~\ref{thm2}. Thus, $\Delta$ is the measure of deviation
between the distributions of $X$ and $Y$, expressed in form of a weighted $L^2$-distance of the corresponding characteristic functions, and this
measure of deviation is
estimated by $T_{n,m}$. As the next result shows, the statistic $T_{n,m}$, when suitably normalized, has a normal limit distribution as $n,m \to \infty$ if $\Delta >0$.
To prove this result, we need the condition
\begin{equation}\label{tsreg}
\lim_{m,n \to \infty} \frac{m}{m+n} = \tau
\end{equation}
for some $\tau \in [0,1]$. In contrast to Theorem~\ref{thm1}  and Theorem~\ref{thm2}, this condition is needed now to assess the asymptotic proportions of
the $X$- sample and the $Y$-sample. In what follows, put
\begin{equation}\label{defk1k2}
K_1(s,t) = \text{Cov}(\text{CS}(s^\top X),\text{CS}(t^\top X)), \quad K_2(s,t) = \text{Cov}(\text{CS}(s^\top Y),\text{CS}(t^\top Y)),
\end{equation}
and let
\begin{equation}\label{defkstar}
K^*(s,t) = \tau K_1(s,t)   +(1-\tau) K_2(s,t), \quad s,t \in M.
\end{equation}
Furthermore, let
\begin{equation}\label{defzt}
z(t) = \varrho(t) - b(t),
\end{equation}
where $\varrho(t)$ and $b(t)$ are given in \eqref{defatbt}.

\bigskip

\begin{thm} \label{thm3} Suppose the standing assumptions on $M$ and $\mu$ hold. If $\Delta >0$, then
\[
\sqrt{\frac{mn}{m+n}} \left(\frac{T_{n,m}}{\frac{mn}{m+n}} - \Delta \right) \verk {\rm{N}}(0,\sigma^2)
\]
under the limiting regime \eqref{tsreg},
where
\[
\sigma^2 = 4 \int_M \int_M K^*(s,t) z(s)z(t) \, \mu({\rm d}s) \mu({\rm d}t).
\]
\end{thm}
\smallskip

{\sc Proof.} The  proof follows the lines of the proof of Theorem 3 of \cite{BEH:2017}. In view of
Theorem~\ref{thm2}, condition (22) of \cite{BEH:2017} holds. Let $A_n$ and $B_m$ as in \eqref{defanbm}, and put $Z_{n,m} := A_n-B_m$.
Furthermore, write $z = z(\cdot)$ for the degenerate random element of  $\HH$, where $z(t)$ is given in \eqref{defzt}, and define
\[
c_{m,n} = \sqrt{\frac{mn}{m+n}}.
\]
Notice that
\[
c_{m,n} \big(Z_{n,m} - z\big) = a_{m,n} \sqrt{n}(A_n-a) - b_{m,n} \sqrt{m}(B_m-b),
\]
where $a_{m,n}$ and $b_{m,n}$ are given in \eqref{defamnbmn}.
By the central limit theorem for $\HH$-valued random elements, we have
 $\sqrt{n}(A_n-a) \verk A$ as $n \to \infty$ and $\sqrt{m}(B_m-b) \verk B$ as $m \to \infty$, where
 $A$ and $B$ are independent centred Gaussian random elements of $\HH$ with covariance kernels $K_1$ and $K_2$, respectively,
 where $K_1$ and $K_2$ are given in \eqref{defk1k2}. In view of \eqref{tsreg}, the continuous mapping theorem yields
 $c_{m,n}(Z_{n,m} - z) \verk Z$, where $Z := \sqrt{\tau} A- \sqrt{1-\tau}B$ is a centred Gaussian random element of $\HH$ having covariance kernel $K^*$
 given in \eqref{defkstar}. Thus, also condition (23) of \cite{BEH:2017} holds, and the proof of Theorem~\ref{thm3} follows in view of
 \[
 c_{m,n} \left(\frac{T_{n,m}}{c^2_{n,m}}- \Delta \right) = 2 \langle c_{m,n}(Z_{n,m} - z),z\rangle + \frac{1}{c_{m,n}} \|c_{m,n}(Z_{n,m}-z)\|^2_\HH.
 \]
Notice that the second summand on the right hand side is $o_\PP(1)$ in view of the tightness of $(c_{m,n}(Z_{n,m}-z))$, and the first summand
converges in distribution to $2\langle Z,z\rangle$, which has the stated normal distribution $\text{N}(0,\sigma^2)$. \bewend
\begin{rem}
Compared to \cite{chen22}, who address the problem of composite hypotheses, the limit results of this section are obtained for simple hypotheses without estimated parameters. However, the results obtained herein hold for artificial sample size $m\neq n$, which is much more general and thus flexible than the case $m=n$ treated by \cite{chen22}. Moreover, our setting is different from that of elliptical distributions on the classical Euclidean space $\mathbb R^d$. In the following section, we suggest a resampling version of the test, and we prove its asymptotic validity.
\end{rem}

\bigskip

\section{Resampling under a simple hypothesis} \label{boot}

Since, under $H_0$, both the finite-sample and the limit distribution of $T_{n,m}$ as $n,m \to \infty$ depend on the unknown underlying
distribution of $X_1$, we use a bootstrap procedure in order to carry out a test that rejects $H_0$ for large values of $T_{n,m}$.
The bootstrap distribution of $T_{n,m}$ is the conditional distribution of $T_{n,m}$ given the pooled sample $X_1,\ldots,X_n,Y_1,\ldots,Y_m$, and a test
of $H_0$ at nominal level $\alpha$ rejects $H_0$ if $T_{n,m}$ exceeds the $(1-\alpha)$-quantile of this bootstrap distribution. Since the bootstrap distribution
is difficult to compute, it is estimated by a Monte Carlo procedure that repeatedly samples from the empirical distribution of the pooled sample.
To be specific, one first computes the observed value $t_{n,m}$ of $T_{n,m}$ based on realizations $x_1,\ldots,x_n,y_1,\ldots,y_m$ of
$X_1,\ldots,X_n,Y_1,\ldots,Y_m$, respectively. In a second step, one generates $b$ independent samples by Monte Carlo simulation. Here, for each $j \in \{1,\ldots,b\}$,
the $j^{{\rm th}}$ sample consists of $x_1(j),\ldots,x_n(j),y_1(j),\ldots,y_m(j)$, where these values have been chosen independently of each other with a uniform distribution over
$\{x_1,\ldots,x_n,y_1,\ldots,y_m\}$. For each $j \in \{1,\ldots,b\}$, one then computes the value $t_{n,m}(j) = T_{n,m}(x_1(j),\ldots,x_n(j),y_1(j),\ldots,y_m(j))$ of the test statistic $T_{n,m}$.
Letting $c_{n,m;1-\alpha}$ denote the $(1-\alpha)$-quantile of the $b$ values $t_{n,m}(1), \ldots, t_{n,m}(b)$, the hypothesis $H_0$ is rejected if
$t_{n,m} > c_{n,m;1-\alpha}$.


To prove that this bootstrap procedure yields a test of $H_0$ of asymptotic level $\alpha$, we use a Hilbert space central limit theorem for
triangular arrays (see \cite{KMM:2000}). This theorem reads as follows.

\begin{thm}\label{thmkundu} Let $\{e_j: j \ge 1\}$ a complete orthonormal basis of the separable Hilbert space $\HH$ with inner product $\langle \cdot, \cdot \rangle $ and norm
$\|\cdot \|_{\mathbb{H}}$. For each $m \ge 1$, let $X_{m1},X_{m2}, \ldots, X_{mm}$ be independent $\HH$-valued random elements such that
$\BE (\langle X_{mj}, e_\ell\rangle ) =0$ and $\BE \|X_{mj}\|_{\mathbb{H}}^2 < \infty$ for each $j \in\{1,\ldots,m\}$ and each $\ell \ge 1$. Put
$S_m = \sum_{j=1}^m X_{mj}$, and let $C_m$ be the covariance operator of $S_m$. Assume that the following conditions hold:
\begin{enumerate}
\item[(i)] $\lim_{m\to \infty} \langle C_me_k,e_\ell\rangle = a_{k\ell}$ (say) exists for each $k \ge 1$ and $\ell \ge 1$,
\item[(ii)] $\lim_{m \to \infty} \sum_{k=1}^\infty \langle C_me_ke_k\rangle = \sum_{k=1}^\infty a_{kk} < \infty$,
\item[(iii)] $\lim_{m\to \infty} L_m(\varepsilon,e_k) = 0$ for each $\varepsilon >0$ and each $k \ge 1$, where
\[
L_m(\varepsilon,h) := \sum_{j=1}^m \BE\big(\langle X_{mj},h \rangle^2 {\bf 1}\big\{ |\langle X_{mj}, h  \rangle | > \varepsilon \big\} \big), \quad h \in \HH .
\]
\end{enumerate}
Then $S_m \verk G$ as $m \to \infty$, where $G$ is a centred random element of $\HH$ with covariance operator $C$ characterized by
$\langle C h,e_\ell \rangle = \sum_{j=1}^\infty \langle h,e_j\rangle a_{j\ell}$ for each $h \in \HH$ and each $\ell \ge 1$.
\end{thm}

\bigskip

To apply Theorem~\ref{thmkundu} in our situation of the Hilbert space $\HH := \L^2(M,{\cal B}(M),\mu)$
 let, in greater generality than considered so far,  $X_1^{m,n}, \ldots, X_n^{m,n},Y_1^{m,n},\ldots, Y_m^{m,n}$ be i.i.d. $M$-valued random vectors with common distribution, and put
\[
\varrho_{m,n}(t) = \BE\big[ \text{CS}\big(t^\top X_1^{m,n}\big)\big], \quad t \in M.
\]
Moreover, let
\begin{eqnarray*}
U_{m,n}(t) & := & \frac{1}{\sqrt{n}} \sum_{j=1}^n \big\{ \text{CS}\big(t^\top X_j^{m,n}\big) - \varrho_{m,n}(t)\big\}, \\
 V_{m,n}(t) & := & \frac{1}{\sqrt{m}} \sum_{k=1}^m \big\{ \text{CS}\big(t^\top Y_k^{m,n}\big) - \varrho_{m,n}(t)\big\}, \quad t \in M
\end{eqnarray*}
(cf. \eqref{defunvm}), and write $U_{m,n} = U_{m,n}(\cdot)$, $V_{m,n} = V_{m,n}(\cdot)$ as well as
\begin{equation}\label{defizmn}
\widetilde{W}_{n,m} := a_{m,n}U_{m,n} - b_{m,n} V_{m,n},
\end{equation}
where $a_{m,n}$ and $b_{m,n}$ are defined in \eqref{defamnbmn}.

\bigskip

\begin{thm}\label{thmbasicsboot} In the setting given above suppose that, under the limiting regime \eqref{tsreg}, we have
$X_1^{m,n} \verk X_\infty$
for some $M$-valued random vector $X_\infty$ with distribution $H_\infty$. Put $\varrho_\infty(t) := \BE[{\rm{CS}}(t^\top X_\infty)]$, $t \in M$, and
$W_\infty(t) := {\rm{CS}}(t^\top X_\infty) - \varrho_\infty(t)$, $t \in M$. Moreover, let $W_\infty := W_\infty(\cdot)$ be the random
element of  $\HH = \L^2(M,{\cal B}(M),\mu)$ with covariance operator $C_\infty$ that is associated with the covariance function
\[
c_\infty(s,t) = \BE \big[{\rm{CS}}(s^\top X_\infty) {\rm{CS}}(t^\top X_\infty) \big] - \varrho_\infty(s) \varrho_\infty(t), \quad s,t \in M,
\]
via
\[
\langle C_\infty g,h\rangle = \int_M \int_M c_\infty (s,t) \, g(s)h(t)\, \mu(\text{d}s) \mu(\text{d}t).
\]
We then have $\widetilde{W}_{n,m} \verk W_\infty$.
\end{thm}

\medskip
{\sc Proof.} The proof is similar to that of Theorem 2 of \cite{BK:2015} and will thus only be sketched.
Notice that
\begin{equation}\label{defxmn1}
X_{m,n,j} := \frac{1}{\sqrt{n}} \left(\text{CS}(\cdot^\top X_j^{m,n})-\rho_{m,n}(\cdot)\right), \quad j \in \{1,\ldots,n\},
\end{equation}
 are i.i.d.
centred random elements of the Hilbert space $\HH = \L^2(M,{\cal B}(M),\mu)$ that, for a fixed complete orthonormal system of $\HH$,  satisfy
$\BE (\langle X_{m,n,j}, e_\ell\rangle ) =0$ and $\BE \|X_{m,n,j}\|_{\mathbb{H}}^2 < \infty$ for each $j \in    \{1,\ldots,n\}$ and each $\ell \ge 1$.
The covariance function of the process $U_{m,n}= U_{m,n}(\cdot)$ is given by
\begin{eqnarray*}
c_{m,n}(s,t) & := & \text{Cov}(U_{m,n}(s),U_{m,n}(t))\\
& = & \BE[ \text{CS}(s^\top X_1^{m,n})\text{CS}(t^\top X_1^{m,n})] - \varrho_{m,n}(s)\varrho_{m,n}(t), \quad s,t \in M,
\end{eqnarray*}
and the covariance operator $C_{m,n}$ (say) of $U_{m,n}$ satisfies
\[
\langle C_{m,n}g,h \rangle = \int_M \int_M c_{m,n}(s,t) \, g(s)h(t) \, \mu(\text{d}s) \mu(\text{d}t), \quad g,h \in \HH.
\]
Since the function $\text{CS}$ is bounded and continuous, $X_1^{m,n} \verk X_\infty$ yields  $c_{m,n}(s,t) \to c_\infty(s,t)$, $s,t \in M$.
By dominated convergence,  we obtain $\langle C_{m,n}g,h\rangle \to \langle C_\infty g,h\rangle$ ($g, \, h \in \HH$), which shows that condition (i)
of Theorem~\ref{thmkundu} holds.  The proof of  condition (ii) of Theorem~\ref{thmkundu} follows the reasoning given on p. 603 of \cite{BK:2015} by replacing $J_0(2\sqrt{t \cdot})$ with $\text{CS}(t^\top \cdot)$, ${\mathcal H}_{m,n}$ with $\varrho_{m,n}$, ${\mathcal H}_\infty$ with $\varrho_\infty$,
$\nu$ with $\mu$, and the region of integration with $M$. To prove condition (iii) of Theorem~\ref{thmkundu}, notice that,   with $X_{m,n,j}$ defined in \eqref{defxmn1},
the fact that $|\text{CS}(\cdot)| \le 2$ and H\"older's inequality give
\[
\big{|} \langle X_{m,n,j}, h \rangle \big{|} \le \frac{2}{\sqrt{n}} (\mu(M))^{1/2}  \|h\|_{\mathbb{H}}^2, \quad h \in \HH.
\]
Consequently,
\[
\sum_{j=1}^n \BE\big(\langle X_{m,n,j},h \rangle^2 {\bf 1}\big\{ |\langle X_{m,n,j}, h  \rangle | > \varepsilon \big\} \big) \le 4 \mu(M) \|h\|_{\mathbb{H}}^2 \PP(|\langle X_{m,n,1}, h  \rangle | > \varepsilon),
\]
and thus also condition (iii) of Theorem~\ref{thmkundu} holds. According to Theorem~\ref{thmkundu}, we have
$U_{m,n} \verk W_\infty$. In the same way, $V_{m,n} \verk \widetilde{W}_\infty$, where, due to the independence of $U_{m,n}$ and $V_{m,n}$,  $\widetilde{W}_\infty$ is an independent
copy of $W_\infty$. In view of \eqref{defizmn} and the continuous mapping theorem, it follows that $\widetilde{W}_{n,m} \verk \sqrt{\tau} W_\infty -   \sqrt{1-\tau} \widetilde{W}_\infty$.
The latter limit has the same distribution as $W_\infty$. \bewend

\bigskip

Notice that the test statistic $T_{n,m}$ figuring in \eqref{reprtmn}, computed on the random variables $X_1^{m,n}, \ldots, X_n^{m,n},Y_1^{m,n},\ldots,Y_n^{m,n}$,  equals
$\|\widetilde{W}_{n,m}\|_{\mathbb{H}}^2 $, where $\widetilde{W}_{n,m}$ is given in \eqref{defizmn}.  From Theorem~\ref{thmbasicsboot}, we thus have the following corollary.

\medskip

\begin{cor} The limit distribution of the test statistic $T_{n,m}$ under the limiting regime \eqref{tsreg} is that of
$\|W_\infty\|{_{\mathbb{H}}}^2$.
\end{cor}

\medskip

Let $F_n$ and $G_m$ denote the empirical distributions of $X_1,\ldots,X_n$ and $Y_1,\ldots,Y_m$, respectively, and write
$
H_{n,m} := \frac{n}{m+n} F_n + \frac{m}{m+n} G_m
$
for the empirical distribution of the pooled sample $X_1,\ldots,X_n,Y_1,\ldots,Y_m$. By the Glivenko-Cantelli theorem, $H_{n,m}$ converges weakly to
$H_\infty := (1-\tau) F + \tau F = F$  with probability one under the limiting regime \eqref{tsreg}, and thus
the bootstrap distribution of $T_{n,m}$ converges almost surely to the distribution of $\|W_\infty\|_{\mathbb{H}}^2$. The latter distribution coincides  with the distribution
of $\|W\|_{\mathbb{H}}^2$, where $W$ is given in Theorem~\ref{thm1}. This shows the asymptotic validity of the bootstrap.

\begin{rem}
The resampling bootstrap procedure applied herein may also be replaced by a permutation procedure. The validity of the exhaustive permutation (that includes all possible permutations) may be directly obtained by observing that, under the null hypothesis $H_0$, the observations $(x_1,...,x_n,y_1,...,y_m)$ are exchangable.  Another potential resampling scheme may be that of weighted bootstrap; see \cite{Alba17}.
\end{rem}

\section{Resampling under a composite hypothesis} \label{boot1}

Analogously to $T_{n,m,w}$, the limit null distribution of $\widehat T_{n,m,w}\linebreak=\widehat T_{n,w}(x _1,...,x _n,\widehat {y }_{1},...,\widehat {y}_{m})$ depends (in a very complicated way) on unknown quantities, and hence
 it cannot be used to compute critical values and actually carry out the test. To this end,
we consider a parametric bootstrap procedure involving the test statistic in \eqref{ts3} computed on the basis of bootstrap observations from $f_0(\cdot;\vartheta)$, where the parameter $\vartheta$ is replaced by estimators.

More precisely, let $\widehat T_{n,m,w, {\rm{obs}}}$ denote the observed value of the test statistic. For given $\alpha \in (0,1)$, write
  $t^*_{n,m,\alpha}$ for the  upper $\alpha$-percentile of the bootstrap distribution of $T_{n,m,w}$. We then define the test
function as
\be \label{psi} \Xi_{n,m}^*=\left\{
         \begin{array}{ll}
           1, & \text{if}\quad T_{n,m,w,{\rm{obs}}} \geq t^*_{n,m,\alpha}, \\
           0, & \text{otherwise}.
         \end{array}
       \right.\ee


 In practice, the  bootstrap distribution of $T_{n,m,w, {\rm{obs}}}$ is approximated as follows:
\begin{enumerate} \itemsep=0pt
\item Generate a bootstrap sample $x _1^*, \ldots, x _n^*$ from $f_0(\cdot;\widehat {\vartheta}_n)$.
 \item Calculate the estimator $\vartheta^*_n=\vartheta_n(x _1^*, \ldots, x _n^*)$
 \item Generate a bootstrap sample $y _{1}^*, \ldots, y _{m}^*$ from $f_0(\cdot;{\vartheta}^*_n)$.
 \item Compute $T_{n,m,w}^*=T_{n,m,w}(x _1^*, \ldots, x _n^*,  y _{1}^*, \ldots, y _{m}^*)$.
\item Repeat steps 1--4 a number of times, say $b$,  and thus  obtain (conditionally on $\widehat{\vartheta}_n$) iid replications of $T^*_{n,,m,w}$, namely $T^*_{n,,m,w,1}, \ldots, T^*_{n,m,w,{\rm{b}}}$.
\end{enumerate}
Then we approximate the  upper $\alpha$--percentile $t^*_{n,m,\alpha}$ in \eqref{psi} of the null distribution of  $T_{n,m,w}$ by the  upper $\alpha$-percentile of the empirical distribution of $T^*_{n,m,w,1}, \ldots, T^*_{n,m,w,{\rm{b}}}$.

Although we provide no asymptotic theory  for the resampling under a composite hypothesis, our simulations show that the above method works well. Nevertheless, it remains an open problem to formally prove that this bootstrap is asymptotically valid.

\section{Simulations} \label{simul}

In this section we provide results of competitive Monte Carlo simulations for the case of both a simple and a composite hypothesis. We throughout restrict the simulation to the spherical setting for the dimension $d=3$ and, for computational feasibility, to the sample size $n=50$. All simulations are performed using the statistical programming language \texttt{R}, see \cite{R:2022}. We implement the test statistic by fixing the measure $\mu(\text{d}t)$ to be the density of the zero-mean spherical stable distribution. Then the test may be computed as in eqn. \eqref{eq:test} with
C$(x)=e^{-\gamma \|x\|^\xi}$, where $(\xi,\gamma) \in (0,2] \times (0,\infty)$ denote tuning parameters which are at our disposal and provide a certain flexibility of the test with respect to power against different alternatives.  Another option, although not yielding a proper measure, is to adopt the approach taken in \cite{SR:2013}, which again results in the test statistic given in \eqref{eq:test} with C$(x)=-\|x\|^\xi, \ \xi\in(0,2)$.


The spherical distributions were generated using the package \texttt{Directional}, see \cite{TASAW:2021}, and the uniformity tests by the package \texttt{sphunif}, see \cite{GV:2020}.

\subsection{Testing the simple hypothesis of uniformity}\label{subsec:unif}
We test the hypothesis
\begin{equation*}
H_0: f(\cdot)\equiv1/|\mathcal{S}^{d-1}|,
\end{equation*}
where $|\mathcal{S}^{d-1}|=2 \pi^{d/2}/\Gamma(d/2)$ is the surface area of the $(d-1)$-sphere. Hence we test whether $f$ is the density of the uniform law $\mathcal{U}(\mathcal{S}^{d-1})$, which is a classical testing problem in directional statistics. For an overview of existing procedures we refer to \cite{GV:2018}.
As competing tests we consider the following procedures:
 \begin{itemize}
 \item The modified Rayleigh test $R_n$, see \cite{MJ:2000}, Section 10.4.1 based on the mean of the directions,

 \item the \cite{EHY:2018} test $NN^J_a$ based on volumes of the $J$ nearest neighbor balls with power $a$,

 \item the Bingham test $B_n$, see \cite{B:1974}, based on the empirical scatter matrix of the sample,

 \item the Sobolev test, see \cite{G:1975}, $G_n$ based on Sobolev norms, and

 \item the test of \cite{CCF:2009} $CA_n$, which is based on random projections that characterize the uniform law.
 \end{itemize}
Empirical critical values for each testing procedure have been obtained by a Monte Carlo simulation study under $H_0$ with 100000 replications.

\begin{figure}[!t]
\centering
\includegraphics[scale=.55]{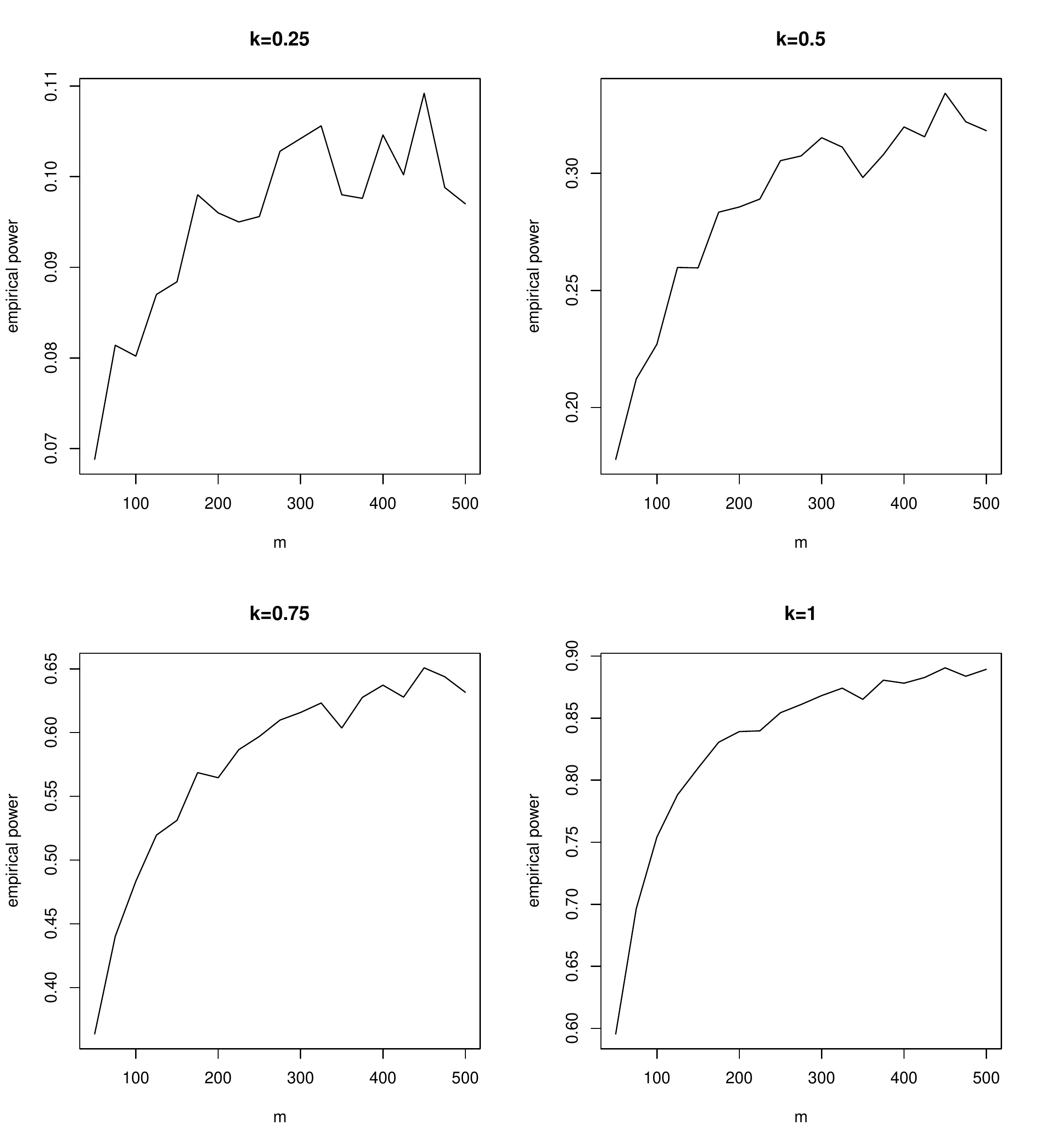}
\caption{Simulated empirical rejection rates for $T^{\{2\}}_{n,1}$ and different values of $10\le m\le 500$ and concentration parameters $\kappa\in\{0.25,0.5,0.75,1\}$ ($n=50$, $\alpha=0.05$, 10000 repetitions) }\label{fig:m}
\end{figure}

We considered the following alternatives to the uniform distribution on $\mathcal{S}^{d-1}$. These alternatives are chosen to simulate different uni-, bi- and trimodal models. For details on the hyperspherical von Mises--Fisher distribution, see Section 9.3 in \cite{MJ:2000}.
\begin{itemize}
    \item The density of the von Mises--Fisher distribution depends on the mean direction $\theta\in\mathcal{S}^{d-1}$ and a concentration parameter $\kappa\ge0$, and it is given by
    \begin{equation*}
        f(x)=\frac{\left( \kappa/2 \right)^{d/2-1}}{\Gamma(d/2) I_{d/2-1}(\kappa)} \exp(\kappa x^\top \theta),\quad x\in\mathcal{S}^{d-1}.
    \end{equation*}
   Here, $I_{d/2-1}$ is the modified Bessel function of the first kind and order $d/2-1$. This class is denoted with $\text{vMF}(\theta, \kappa)$.

    \item We simulate a mixture of two von Mises--Fisher distributions with different mean directions by the following procedure. Simulate $U\sim\mathcal{U}(0,1)$ and, independently, $Y_i \sim \text{vMF}(\theta_i, \kappa_i)$, $i\in\{1,2\}$, with corresponding location and concentration parameters, and choose $p \in (0,1)$. Then we generate a member $X$ of the random sample according to
    \begin{equation*}
        X = Y_1 \mathbf{1}{\{ U < p \} } + Y_2 \mathbf{1}{\{ U \geq p \} }.
    \end{equation*}
    We denote this alternative class with $\text{MMF}((p,1-p), (\theta_1, \theta_2), (\kappa_1, \kappa_2))$.

    \item In a similar manner as for the mixture of two vMF distributions, we simulate a mixture of three  von Mises--Fisher distributions with different centers, by additionally simulating independently a third random vector $Y_3\sim \text{vMF}(\theta_3, \kappa_3)$ and generating the member $X$ by
    \begin{equation*}
        X = Y_1 \mathbf{1}{\{ U < p \} } + Y_2 \mathbf{1}{\{ p \leq U < 2p \} } + Y_3 \mathbf{1}{\{ U \geq 2p \} }.
    \end{equation*}
    We denote this class with $\text{MMF}((p,p,1-2p), (\theta_1, \theta_2, \theta_3), (\kappa_1, \kappa_2, \kappa_3))$.
\end{itemize}
In each of the alternatives, we put $\mu_1=(1,0,\ldots,0),\, \mathtt{1}=(1,\ldots,1)/\sqrt{d}$, and \linebreak$\mu_2=(-1,1,\ldots,1)/\sqrt{d}$.
Here and in the following, $T_{n,\gamma}^{\{\xi\}}$ stands for the test in eqn. \eqref{eq:test} with C$(x)=e^{-\gamma \|x\|^\xi}$, where $(\xi,\gamma) \in (0,2] \times (0,\infty)$ as well as $T^{SR}_{n,a}$  for the test with C$(x)=-\|x\|^a$. The result of the simulation is displayed in Tables \ref{tab:unif1} and \ref{tab:unif12} for the choice $\xi=2$ and the stated competitors. As can be seen, the suggested tests perform well in comparison, although they are never the best performing procedures. This behavior might be explained by the approximation of the true characteristic function under the null hypothesis. To investigate the impact of the sample size $m$ of the simulated data set $Y_1,\ldots,Y_m$, we simulated the empirical power of the test for four vMF distributions and for different values of $m$, see Figure \ref{fig:m}. Clearly, the choice of $m$ has an impact on the estimation, and larger values of $m$  are desirable, but increasing $m$ leads to longer computation time. Table \ref{tab:unif2} exhibits the impact of the weighting measure $\mu$ and hence of the choice of the function C$(\cdot)$. In terms of power for the  uni- and bimodal alternatives considered, the choice of C$(\cdot)$ has nearly no influence on the empirical power, with the exception of the MMF$((0.5,0.5),(-\mu_1,\mu_1),(2,2))$ alternative, where $T_{n,\gamma}^{\{1\}}$, $T_{n,\gamma}^{\{1.5\}}$ and $T_{n,\gamma}^{\{2\}}$ outperform the
$T^{SR}_{n,a}$-procedures for some values of the tuning parameter $\gamma$.

\begin{table}[t]
\centering
\footnotesize
\caption{Empirical rejection rates for testing uniformity for the test $T^{\{2\}}_{n,\gamma}$ and competitors ($n=50$, $m=500$, $\alpha=0.05$, 10000 replications)} \label{tab:unif1}
	\setlength{\tabcolsep}{.8mm}
\begin{tabular}{l|rrrrrrrrrrr}
  Alternative & $T^{\{2\}}_{n,5}$ & $T^{\{2\}}_{n,2}$ & $T^{\{2\}}_{n,1}$ & $T^{\{2\}}_{n,.5}$ & $T^{\{2\}}_{n,.25}$ & $T^{\{2\}}_{n,.17}$ &  $R_n$ & $NN^{15}_{.5}$  & $B_n$ & $G_n$ & $CA$ \\
  \hline
vMF$(\mu_1,0)=\mathcal{U}(\mathcal{S}^2)$ & 5 & 5 & 5 & 5 & 5 & 4 & 5 & 4 & 5 & 5 & 5 \\
vMF$(\mu_1,.25)$ & 8 & 11 & 11 & 11 & 11 & 11 & 12 & 8 & 5 & 12 & 11 \\
vMF$(\mu_1,.5)$ & 20 & 31 & 32 & 32 & 33 & 33 & 37 & 21 & 5 & 36 & 32 \\
vMF$(\mu_1,.75)$ & 44 & 63 & 64 & 66 & 67 & 66 & 72 & 49 & 7 & 71 & 64 \\
vMF$(\mu_1,1)$ & 71 & 88 & 89 & 89 & 90 & 89 & 93 & 77 & 9 & 92 & 89 \\
MMF$((.25,.75),(\mu_1,\mu_1),(0,2))$ & 95 & 99 & 99 & 99 & 99 & 99 & 99 & 96 & 39 & 99 & 99 \\
MMF$((.25,.75),(-\mu_1,\mu_1),(2,2))$ & 78 & 83 & 81 & 80 & 77 & 75 & 77 & 78 & 64 & 85 & 81 \\
MMF$((.5,.5),(\mu_1,\mu_1),(0,2))$ & 60 & 75 & 76 & 76 & 76 & 75 & 80 & 61 & 18 & 81 & 75 \\
MMF$((.5,.5),(-\mu_1,\mu_1),(2,2))$ & 37 & 23 & 16 & 12 & 8 & 7 & 6 & 48 & 63 & 18 & 13 \\
MMF$((.75,.25),(\mu_1,\mu_1),(0,2))$ & 17 & 23 & 23 & 24 & 24 & 23 & 26 & 16 & 8 & 26 & 23 \\
MMF$((.75,.25),(-\mu_1,\mu_1),(2,2))$ & 76 & 83 & 82 & 80 & 76 & 75 & 77 & 78 & 63 & 85 & 82 \\
MMF$((.25,.75),(-\mu_1,\mu_1),(5,0))$ & 48 & 54 & 55 & 52 & 50 & 49 & 52 & 48 & 32 & 57 & 51 \\
MMF$((.25,.75),(-\mu_1,\mu_1),(5,1))$ & 35 & 21 & 15 & 12 & 9 & 7 & 7 & 44 & 53 & 17 & 13 \\
MMF$((.25,.75),(-\mu_1,\mu_1),(5,2))$ & 88 & 83 & 76 & 70 & 59 & 53 & 51 & 93 & 93 & 81 & 74 \\
MMF$((.25,.75),(-\mu_1,\mu_1),(5,3))$ & 100 & 100 & 99 & 98 & 91 & 88 & 83 & 100 & 100 & 99 & 98 \\
MMF$((.25,.75),(-\mu_1,\mu_1),(5,4))$ & 100 & 100 & 100 & 100 & 98 & 97 & 94 & 100 & 100 & 100 & 100 \\
MMF$((.25,.75),(-\mu_1,\mu_1),(0,3))$ & 100 & 100 & 100 & 100 & 100 & 100 & 100 & 100 & 86 & 100 & 100 \\
MMF$((.25,.75),(-\mu_1,\mu_1),(1,3))$ & 99 & 100 & 100 & 100 & 99 & 99 & 99 & 99 & 90 & 100 & 100 \\
MMF$((.25,.75),(-\mu_1,\mu_1),(2,3))$ & 98 & 99 & 98 & 98 & 97 & 96 & 96 & 99 & 96 & 99 & 98 \\
MMF$((.25,.75),(-\mu_1,\mu_1),(3,3))$ & 99 & 99 & 98 & 97 & 94 & 92 & 91 & 99 & 99 & 99 & 98 \\
MMF$((.25,.75),(-\mu_1,\mu_1),(4,3))$ & 100 & 100 & 100 & 100 & 100 & 100 & 100 & 100 & 76 & 100 & 100 \\
MMF$((.5,.5),(-\mu_1,\mu_1),(5,0))$ & 98 & 99 & 99 & 99 & 99 & 98 & 99 & 98 & 89 & 99 & 99 \\
MMF$((.5,.5),(-\mu_1,\mu_1),(5,1))$ & 93 & 89 & 85 & 81 & 72 & 68 & 66 & 94 & 95 & 88 & 84 \\
MMF$((.5,.5),(-\mu_1,\mu_1),(5,2))$ & 96 & 90 & 81 & 68 & 42 & 34 & 25 & 99 & 100 & 85 & 71 \\
MMF$((.5,.5),(-\mu_1,\mu_1),(5,3))$ & 99 & 98 & 93 & 82 & 36 & 22 & 11 & 100 & 100 & 96 & 76 \\
MMF$((.5,.5),(-\mu_1,\mu_1),(5,4))$ & 100 & 100 & 99 & 95 & 46 & 23 & 8 & 100 & 100 & 100 & 89 \\
MMF$((.5,.5),(-\mu_1,\mu_1),(0,3))$ & 85 & 93 & 93 & 93 & 93 & 92 & 94 & 86 & 48 & 95 & 92 \\
MMF$((.5,.5),(-\mu_1,\mu_1),(1,3))$ & 60 & 58 & 53 & 51 & 45 & 42 & 43 & 63 & 64 & 58 & 52 \\
MMF$((.5,.5),(-\mu_1,\mu_1),(2,3))$ & 69 & 53 & 38 & 29 & 17 & 14 & 11 & 80 & 90 & 44 & 31 \\
MMF$((.5,.5),(-\mu_1,\mu_1),(3,3))$ & 89 & 77 & 59 & 42 & 17 & 11 & 6 & 96 & 99 & 65 & 41 \\
MMF$((.5,.5),(-\mu_1,\mu_1),(4,3))$ & 100 & 100 & 100 & 100 & 100 & 100 & 100 & 100 & 71 & 100 & 100
\end{tabular}
\end{table}

\begin{landscape}
\begin{table}[t]
\centering
\footnotesize
\caption{Empirical rejection rates for testing uniformity for the test $T^{\{2\}}_{n,\gamma}$ and competitors ($n=50$, $m=500$, $\alpha=0.05$, 10000 replications)} \label{tab:unif12}
	\setlength{\tabcolsep}{.8mm}
\begin{tabular}{l|rrrrrrrrrrr}
  Alternative & $T^{\{2\}}_{n,5}$ & $T^{\{2\}}_{n,2}$ & $T^{\{2\}}_{n,1}$ & $T^{\{2\}}_{n,.5}$ & $T^{\{2\}}_{n,.25}$ & $T^{\{2\}}_{n,.17}$ &  $R_n$ & $NN^{15}_{.5}$  & $B_n$ & $G_n$ & $CA$ \\
  \hline
MMF$((.25,.25,.5),(-\mathtt{1},\mu_2,\mu_1),(2,2,2))$ & 23 & 22 & 21 & 20 & 19 & 18 & 19 & 25 & 23 & 24 & 20 \\
MMF$((.25,.25,.5),(-\mathtt{1},\mu_2,\mu_1),(2,3,2))$ & 25 & 21 & 18 & 18 & 16 & 14 & 15 & 29 & 29 & 21 & 17 \\
MMF$((.25,.25,.5),(-\mathtt{1},\mu_2,\mu_1),(2,4,2))$ & 32 & 24 & 20 & 18 & 16 & 14 & 15 & 38 & 38 & 23 & 17 \\
MMF$((.25,.25,.5),(-\mathtt{1},\mu_2,\mu_1),(2,2,1))$ & 7 & 6 & 6 & 5 & 5 & 5 & 5 & 8 & 10 & 6 & 5 \\
MMF$((.25,.25,.5),(-\mathtt{1},\mu_2,\mu_1),(2,2,3))$ & 56 & 55 & 51 & 48 & 44 & 42 & 43 & 58 & 54 & 56 & 49 \\
MMF$((.25,.25,.5),(-\mathtt{1},\mu_2,\mu_1),(2,2,4))$ & 79 & 76 & 71 & 68 & 61 & 58 & 58 & 82 & 79 & 76 & 69 \\
MMF$((.25,.25,.5),(-\mathtt{1},\mu_2,\mu_1),(.5,2,1))$ & 11 & 13 & 13 & 14 & 13 & 12 & 14 & 12 & 9 & 15 & 13 \\
MMF$((.25,.25,.5),(-\mathtt{1},\mu_2,\mu_1),(.25,1,2))$ & 44 & 57 & 56 & 57 & 56 & 54 & 60 & 45 & 19 & 62 & 55 \\
MMF$((.25,.25,.5),(-\mathtt{1},\mu_2,\mu_1),(2,.25,1))$ & 14 & 17 & 18 & 18 & 18 & 17 & 20 & 16 & 9 & 20 & 17 \\
MMF$((.1,.1,.8),(-\mathtt{1},\mu_2,\mu_1),(2,2,2))$ & 91 & 97 & 97 & 97 & 97 & 97 & 98 & 93 & 44 & 98 & 97 \\
MMF$((.1,.1,.8),(-\mathtt{1},\mu_2,\mu_1),(2,3,2))$ & 91 & 96 & 97 & 97 & 97 & 96 & 97 & 92 & 46 & 98 & 97 \\
MMF$((.1,.1,.8),(-\mathtt{1},\mu_2,\mu_1),(2,4,2))$ & 90 & 96 & 96 & 96 & 95 & 95 & 97 & 91 & 47 & 98 & 95 \\
MMF$((.1,.1,.8),(-\mathtt{1},\mu_2,\mu_1),(2,2,1))$ & 29 & 41 & 42 & 43 & 43 & 42 & 48 & 29 & 8 & 48 & 41 \\
MMF$((.1,.1,.8),(-\mathtt{1},\mu_2,\mu_1),(2,2,3))$ & 100 & 100 & 100 & 100 & 100 & 100 & 100 & 100 & 91 & 100 & 100 \\
MMF$((.1,.1,.8),(-\mathtt{1},\mu_2,\mu_1),(2,2,4))$ & 100 & 100 & 100 & 100 & 100 & 100 & 100 & 100 & 100 & 100 & 100 \\
MMF$((.1,.1,.8),(-\mathtt{1},\mu_2,\mu_1),(.5,2,1))$ & 37 & 52 & 53 & 54 & 55 & 53 & 59 & 39 & 8 & 59 & 52 \\
MMF$((.1,.1,.8),(-\mathtt{1},\mu_2,\mu_1),(.25,1,2))$ & 96 & 99 & 99 & 99 & 99 & 99 & 100 & 97 & 43 & 100 & 99 \\
MMF$((.1,.1,.8),(-\mathtt{1},\mu_2,\mu_1),(2,.25,1))$ & 36 & 53 & 55 & 55 & 55 & 54 & 60 & 40 & 9 & 60 & 53
\end{tabular}
\end{table}
\end{landscape}

\begin{landscape}
\begin{table}[t]
\centering
\footnotesize
\setlength{\tabcolsep}{.7mm}
\caption{Empirical rejection rates for testing uniformity for the test $T^{\{\xi\}}_{n,\gamma}$ for $\xi=1,1.5$ as well as $T^{SR}_{n,a}$ ($n=50$, $m=500$, $\alpha=0.05$, 10000 replications)} \label{tab:unif2}
\begin{tabular}{l|rrrrrrrrrrrrrrrrrrr}
Alternative & $T^{\{1\}}_{n,.1}$ & $T^{\{1\}}_{n,.5}$ & $T^{\{1\}}_{n,.75}$ & $T^{\{1\}}_{n,1}$ & $T^{\{1\}}_{n,2}$ & $T^{\{1\}}_{n,3}$ & $T^{\{1.5\}}_{n,.1}$ & $T^{\{1.5\}}_{n,.5}$ & $T^{\{1.5\}}_{n,.75}$ & $T^{\{1.5\}}_{n,1}$ & $T^{\{1.5\}}_{n,2}$ & $T^{\{1.5\}}_{n,3}$ & $T^{SR}_{n,.1}$ & $T^{SR}_{n,.5}$ & $T^{SR}_{n,.75}$ & $T^{SR}_{n,1}$ & $T^{SR}_{n,2}$ & $T^{SR}_{n,3}$\\
  \hline
vMF$(\mu_1,0)=\mathcal{U}(\mathcal{S}^2)$ & 5 & 5 & 5 & 5 & 6 & 5 & 4 & 4 & 5 & 5 & 6 & 5 & 5 & 4 & 5 & 6 & 6 & 5\\
vMF$(\mu_1,.25)$ & 11 & 11 & 12 & 11 & 11 & 10 & 12 & 12 & 13 & 11 & 11 & 10 & 11 & 11 & 12 & 12 & 12 & 12\\
vMF$(\mu_1,.5)$ & 33 & 33 & 32 & 30 & 29 & 25 & 33 & 33 & 33 & 32 & 28 & 26 & 31 & 32 & 33 & 33 & 34 & 34\\
vMF$(\mu_1,.75)$ & 65 & 64 & 64 & 61 & 59 & 52 & 66 & 66 & 65 & 63 & 58 & 53 & 63 & 64 & 66 & 65 & 67 & 66\\
vMF$(\mu_1,1)$ & 90 & 89 & 89 & 88 & 85 & 80 & 90 & 89 & 90 & 88 & 84 & 79 & 87 & 88 & 89 & 89 & 90 & 90\\
MMF$((.25,.75),(\mu_1,\mu_1),(0,2))$  & 99 & 99 & 99 & 99 & 98 & 97 & 99 & 99 & 99 & 99 & 98 & 97 & 99 & 99 & 99 & 99 & 99 & 99\\
MMF$((.25,.75),(-\mu_1,\mu_1),(2,2))$ & 79 & 80 & 82 & 82 & 83 & 81 & 77 & 80 & 82 & 82 & 84 & 82 & 82 & 82 & 81 & 80 & 79 & 76\\
MMF$((.5,.5),(\mu_1,\mu_1),(0,2))$ & 77 & 76 & 76 & 75 & 73 & 67 & 77 & 77 & 77 & 76 & 74 & 68 & 75 & 76 & 76 & 77 & 77 & 76\\
MMF$((.5,.5),(-\mu_1,\mu_1),(2,2))$ & 10 & 13 & 18 & 21 & 33 & 37 & 7 & 12 & 19 & 22 & 36 & 40 & 20 & 17 & 15 & 13 & 11 & 7\\
MMF$((.75,.25),(\mu_1,\mu_1),(0,2))$ & 25 & 25 & 25 & 24 & 23 & 21 & 24 & 25 & 25 & 24 & 23 & 21 & 24 & 25 & 26 & 25 & 25 & 25\\
MMF$((.75,.25),(-\mu_1,\mu_1),(2,2))$ & 79 & 80 & 82 & 83 & 84 & 81 & 76 & 80 & 83 & 83 & 85 & 83 & 81 & 82 & 83 & 81 & 81 & 78\\
MMF$((.25,.75),(-\mu_1,\mu_1),(5,0))$ & 51 & 52 & 55 & 55 & 56 & 55 & 48 & 52 & 55 & 55 & 56 & 54 & 52 & 53 & 54 & 53 & 52 & 49\\
MMF$((.25,.75),(-\mu_1,\mu_1),(5,1))$ & 12 & 15 & 18 & 20 & 32 & 35 & 8 & 13 & 17 & 20 & 33 & 36 & 20 & 17 & 15 & 13 & 12 & 8\\
MMF$((.25,.75),(-\mu_1,\mu_1),(5,2))$ & 66 & 74 & 79 & 81 & 88 & 89 & 58 & 71 & 78 & 81 & 89 & 90 & 81 & 79 & 76 & 70 & 65 & 55\\
MMF$((.25,.75),(-\mu_1,\mu_1),(5,3))$ & 96 & 98 & 99 & 99 & 100 & 100 & 91 & 97 & 99 & 99 & 100 & 100 & 99 & 99 & 98 & 97 & 95 & 88\\
MMF$((.25,.75),(-\mu_1,\mu_1),(5,4))$ & 100 & 100 & 100 & 100 & 100 & 100 & 98 & 100 & 100 & 100 & 100 & 100 & 100 & 100 & 100 & 100 & 100 & 97\\
MMF$((.25,.75),(-\mu_1,\mu_1),(0,3))$ & 100 & 100 & 100 & 100 & 100 & 100 & 100 & 100 & 100 & 100 & 100 & 100 & 100 & 100 & 100 & 100 & 100 & 100\\
MMF$((.25,.75),(-\mu_1,\mu_1),(1,3))$ & 99 & 99 & 99 & 99 & 99 & 99 & 99 & 99 & 100 & 100 & 99 & 99 & 100 & 99 & 100 & 99 & 100 & 9\\
MMF$((.25,.75),(-\mu_1,\mu_1),(2,3))$ & 98 & 98 & 99 & 99 & 99 & 99 & 97 & 98 & 99 & 99 & 99 & 99 & 99 & 99 & 99 & 98 & 98 & 97\\
MMF$((.25,.75),(-\mu_1,\mu_1),(3,3))$ & 96 & 97 & 98 & 98 & 99 & 99 & 94 & 98 & 99 & 99 & 99 & 99 & 99 & 99 & 98 & 97 & 96 & 93\\
MMF$((.25,.75),(-\mu_1,\mu_1),(4,3))$ & 100 & 100 & 100 & 100 & 100 & 100  & 100 & 100 & 100 & 100 & 100 & 100 & 100 & 100 & 100 & 100 & 100 & 100\\
\end{tabular}
\end{table}
\end{landscape}

\subsection{Testing the fit to the von Mises--Fisher distribution}\label{subsec:vMF}
For the case of a composite hypothesis, we consider the hypothesis that the underlying density belongs to the family of von Mises--Fisher distributions vMF$(\kappa,\theta)$, i.e., we test the hypothesis
\begin{equation}\label{eq:H0compvMF}
  H_0:   f(\cdot)=\frac{\left( \kappa/2 \right)^{d/2-1}}{\Gamma(d/2) I_{d/2-1}(\kappa)} \exp(\kappa \cdot^\top \theta),\quad \mbox{for some } \kappa\ge0\mbox{ and } \theta\in\mathcal{S}^{d-1},
\end{equation}
against general alternatives. 
The main difference to subsection \ref{subsec:unif} is that we consider a test to a \textit{family} of distributions, where the parameters are unknown and hence have to be estimated. To test the hypothesis we chose $T_{n,\gamma}^{\{2\}}$ in eqn. \eqref{eq:test} with C$(x)=e^{-\gamma \|x\|^2}$, for different values of the tuning parameter $\gamma$, and we implemented the parametric bootstrap procedure from Section \ref{boot1}. To approximate the unknown parameters we calculated the maximum likelihood estimates for $\kappa$ and $\theta$ as proposed in Section 10.3.1 of \cite{MJ:2000}. As far as we know, testing  composite hypotheses for spherical or hyperspherical distributions with estimated parameters has not been considered before in the literature. As alternative models we chose the same distributions as described in Subsection \ref{subsec:unif}.

In view of the extensive computation time due to the parametric bootstrap procedure, we considered the simulation setting $n=50$, $m=200$, a sample size of $500$ in the bootstrap algorithm, and $5000$ Monte Carlo replications. Throughout the study, we fixed the significance level to 0.05. The results are reported in Table \ref{tab:comp_Hyp}. Notably, the novel test maintains the nominal significance level very closely, and its power with respect to bimodal alternatives increases the more these two modes are pronounced.

\begin{table}[t]
\centering
\caption{Empirical rejection rates for testing the fit to a von Mises--Fisher distribution ($n=50$, $m=200$, $\alpha=0.05$, 500 bootstrap sample size, 5000 replications)  }\label{tab:comp_Hyp}
\begin{tabular}{ l|rrrrrr }
Alternative &  $T_{n,.1}^{\{2\}}$ & $T_{n,.5}^{\{2\}}$ & $T_{n,.75}^{\{2\}}$ & $T_{n,1}^{\{2\}}$ & $T_{n,2}^{\{2\}}$ & $T_{n,3}^{\{2\}}$  \\\hline
vMF$(\mu_1,0)=\mathcal{U}(\mathcal{S}^2)$ & 5 & 6 & 6 & 6 & 6 & 6 \\
vMF$(\mu_1,.25)$ & 5 & 5 & 6 & 6 & 6 & 6 \\
vMF$(\mu_1,.5)$  & 5 & 5 & 6 & 5 & 5 & 6 \\
vMF$(\mu_1,.75)$  & 5 & 5 & 6 & 5 & 6 & 6 \\
vMF$(\mu_1,1)$  & 5 & 5 & 6 & 6 & 5 & 6 \\
MMF$((.25,.75),(\mu_1,\mu_1),(0,2))$ & 6 & 7 & 7 & 7 & 7 & 6 \\
MMF$((.25,.75),(-\mu_1,\mu_1),(2,2))$ & 21 & 31 & 30 & 27 & 19 & 13 \\
MMF$((.5,.5),(\mu_1,\mu_1),(0,2))$  & 6 & 7 & 8 & 8 & 7 & 6 \\
MMF$((.5,.5),(-\mu_1,\mu_1),(2,2))$ & 36 & 49 & 47 & 45 & 30 & 20 \\
MMF$((.75,.25),(\mu_1,\mu_1),(0,2))$ & 5 & 7 & 7 & 7 & 6 & 6 \\
MMF$((.75,.25),(-\mu_1,\mu_1),(2,2))$  & 21 & 31 & 31 & 29 & 20 & 14 \\
MMF$((.25,.75),(-\mu_1,\mu_1),(5,0))$ & 16 & 18 & 17 & 16 & 12 & 9 \\
MMF$((.25,.75),(-\mu_1,\mu_1),(5,1))$ & 36 & 42 & 41 & 38 & 25 & 16 \\
MMF$((.25,.75),(-\mu_1,\mu_1),(5,2))$ & 59 & 73 & 72 & 69 & 54 & 37 \\
MMF$((.25,.75),(-\mu_1,\mu_1),(5,3))$ & 80 & 92 & 92 & 91 & 82 & 65 \\
MMF$((.25,.75),(-\mu_1,\mu_1),(5,4))$ & 92 & 98 & 98 & 98 & 94 & 87 \\
MMF$((.25,.75),(-\mu_1,\mu_1),(0,3))$ & 10 & 13 & 13 & 12 & 9 & 8 \\
MMF$((.25,.75),(-\mu_1,\mu_1),(1,3))$ & 23 & 32 & 31 & 29 & 20 & 15 \\
MMF$((.25,.75),(-\mu_1,\mu_1),(2,3))$& 41 & 59 & 57 & 55 & 39 & 26 \\
MMF$((.25,.75),(-\mu_1,\mu_1),(3,3))$  & 60 & 78 & 78 & 75 & 59 & 44 \\
MMF$((.25,.75),(-\mu_1,\mu_1),(4,3))$ & 25 & 34 & 32 & 30 & 19 & 13
\end{tabular}
\end{table}

\subsection{Testing the fit to the angular central Gaussian distribution}
In this subsection, we consider testing the fit to an angular central Gaussian model, i.e., we test the hypothesis
\begin{equation}\label{eq:H0compacG}
  H_0:   f(\cdot)=|\Sigma|^{-1/2}(\cdot^\top \Sigma^{-1} \cdot)^{-d/2} ,\quad \mbox{for some } \Sigma.
\end{equation}
Here, $\Sigma$ is a symmetric positive definite ($d \times d$)-parameter matrix, which is identifiable up to multiplication by a positive scalar. For information regarding this model, see \cite{MJ:2000}, Section 9.4.4, and for a numerical procedure to approximate the maximum likelihood estimator of the unknown parameter matrix $\Sigma$, see \cite{T:1987}. To the best of our knowledge, testing the fit to the angular central Gaussian family has not been considered in the literature.

The simulation parameters match the ones of Subsection \ref{subsec:vMF}. In complete analogy, we considered $T_{n,\gamma}^{\{2\}}$ in eqn. \eqref{eq:test} with C$(x)=e^{-\gamma \|x\|^2}$, for different values of the tuning parameter $\gamma$, and we implemented the parametric bootstrap procedure from Section \ref{boot1}. To simulate different models under the null hypothesis, we generated a realisation of a random matrix
\begin{equation*}
    A=\left(\begin{array}{ccc}
    -0.846 & -0.531 & -0.779 \\
   0.609 & -0.096 & 0.761 \\
   0.851 & 0.133 & -0.666
   \end{array}\right)
\end{equation*}
and computed the covariance matrices $\Sigma_\ell=(A^\ell)^\top A^\ell$, $\ell\in\{1,2,3,4\}$, where the power matrix $A^\ell=(a_{ij}^\ell)_{i,j=1,2,3}$ is defined component-by-component, and $\Sigma_0 = \text{diag}(1,2,3)$. The corresponding alternatives are denoted by ACG$_\ell$, $\ell=0,\ldots,4$.  Results are presented in Table \ref{tab:comp_Hyp_ACG}. In this case the bootstrap testing procedure controls the type I error, while performing well for most of the alternatives considered.

\begin{table}[ht]
\centering
\caption{Empirical rejection rates for testing the fit to an angular central Gaussian distribution ($n=50$, $m=200$, $\alpha=0.05$, 500 bootstrap sample size, 5000 replications)  }\label{tab:comp_Hyp_ACG}
\begin{tabular}{l|rrrr}
Alternative &  $T_{n,.1}^{\{2\}}$ & $T_{n,.5}^{\{2\}}$ & $T_{n,1}^{\{2\}}$ & $T_{n,3}^{\{2\}}$  \\\hline
ACG$_0$ & 4 & 4 & 4 & 3 \\
ACG$_1$ & 5 & 4 & 4 & 4 \\
ACG$_2$ & 4 & 6 & 4 & 5 \\
ACG$_3$ & 4 & 3 & 5 & 5 \\
ACG$_4$ & 6 & 6 & 5 & 4 \\
vMF$(\mu_1,0)=\mathcal{U}(\mathcal{S}^2)$ & 3 & 3 & 1 & 1 \\
vMF$(\mu_1,.25)$ & 12 & 15 & 14 & 12 \\
vMF$(\mu_1,.5)$  & 23 & 34 & 26 & 18 \\
vMF$(\mu_1,.75)$ & 45 & 74 & 60 & 50 \\
vMF$(\mu_1,1)$ & 74 & 91 & 91 & 79 \\
MMF$((.25,.75),(\mu_1,\mu_1),(0,2))$ & 95 & 100 & 100 & 95 \\
MMF$((.25,.75),(-\mu_1,\mu_1),(2,2))$ & 41 & 69 & 66 & 47 \\
MMF$((.5,.5),(\mu_1,\mu_1),(0,2))$ & 60 & 83 & 72 & 67 \\
MMF$((.5,.5),(-\mu_1,\mu_1),(2,2))$ & 4 & 3 & 3 & 1 \\
MMF$((.75,.25),(\mu_1,\mu_1),(0,2))$ & 6 & 17 & 11 & 9 \\
MMF$((.75,.25),(-\mu_1,\mu_1),(2,2))$ & 54 & 59 & 58 & 76 \\
MMF$((.25,.75),(-\mu_1,\mu_1),(5,0))$ & 25 & 23 & 41 & 46 \\
MMF$((.25,.75),(-\mu_1,\mu_1),(5,1))$ & 2 & 2 & 2 & 0 \\
MMF$((.25,.75),(-\mu_1,\mu_1),(5,2))$& 28 & 44 & 39 & 28 \\
MMF$((.25,.75),(-\mu_1,\mu_1),(5,3))$ & 46 & 62 & 63 & 49 \\
MMF$((.25,.75),(-\mu_1,\mu_1),(5,4))$  & 76 & 91 & 84 & 80 \\
MMF$((.25,.75),(-\mu_1,\mu_1),(0,3))$ & 98 & 100 & 99 & 99 \\
MMF$((.25,.75),(-\mu_1,\mu_1),(1,3))$ & 95 & 99 & 98 & 98 \\
MMF$((.25,.75),(-\mu_1,\mu_1),(2,3))$ & 83 & 94 & 92 & 86 \\
MMF$((.25,.75),(-\mu_1,\mu_1),(3,3))$& 69 & 86 & 81 & 74 \\
MMF$((.25,.75),(-\mu_1,\mu_1),(4,3))$ & 100 & 100 & 100 & 100 
\end{tabular}
\end{table}

\section{Real data} \label{real}

\begin{figure}[t]
\centering
\includegraphics[scale=0.25]{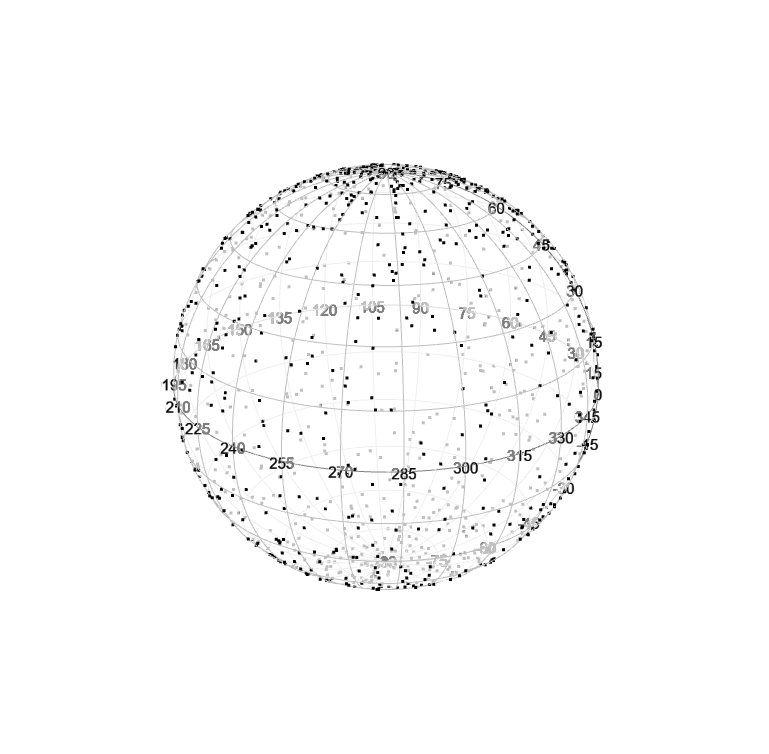}\includegraphics[scale=0.35]{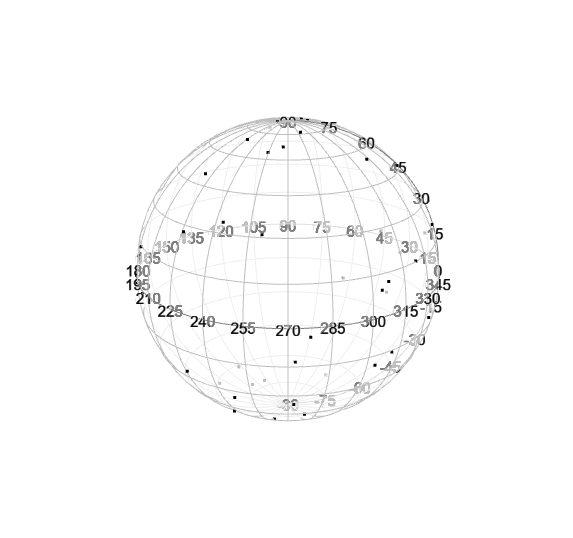}
\caption{Full data set ($n=1137$) of archeomagnetic directions (left) and subsample ($n=50$) consisting of directions the same age 1250 (right) .}\label{fig:pal1}
\end{figure}

We revisit the paleomagnetic data in \cite{SW:2019}, which is an example of spherical data.  Paleomagnetic data consist of observations on the direction of magnetism in either rocks, sediment, or in archeological specimens. These data are measured at various geological points in time and spatial locations. The directions are usually measured as declination and inclination angles based on strike and dip coordinates, see \cite{SW:2019} and the references therein for more information. The data considered are taken from the GEOMAGIA50.v3 database, see \cite{Betal:2015}. For simplicity, we analyse the data provided in the supplementary material of \cite{SW:2019}. The full data set consists of $n=1137$ entries (variables are \texttt{age}, \texttt{dec}, \texttt{inc}, \texttt{lat}, and \texttt{lon}) collected at a single spatial location, which is the Eifel maars (EIF) lakes in Germany with relocated nearby data, for details see \cite{SW:2019}. The analysed directions are given by the variables declination $D$ (\texttt{dec} defined on $[0^\circ,360^\circ]$) and inclination $I$ (\texttt{inc} defined on $[-90^\circ,90^\circ]$). They are converted to Cartesian coordinates by $x_1=\sin(I)$, $x_2=\cos(I)\cos(D)$, and $x_3=\cos(I)\sin(D)$, ensuring $x=(x_1,x_2,x_3)\in\mathcal{S}^2$. For a plot of the data, see Figure \ref{fig:pal1} (left).

We test the composite hypothesis \eqref{eq:H0compvMF} of a von Mises--Fisher distribution by $T_{n,\gamma}^{\{2\}}$ in eqn. \eqref{eq:test} with C$(x)=e^{-\gamma \|x\|^2}$, for different values of the tuning parameter $\gamma$, and we fix $m=500$ with a bootstrap sample size $b=1000$. The bootstrap
$p$-values are reported in Table \ref{tab:pal1}. With the exception of the tuning parameter $\gamma=0.5$, the $p$-values indicate that we are not able to reject the hypothesis of fit of an underlying Mises--Fisher distribution at any level. For their analysis, the authors in \cite{SW:2019} consider  rocks of age 1250 and hence determine a subset of the data of sample size $n=50$, for a plot see Figure \ref{fig:pal1} (right). They propose to use a new spherical model, namely a distribution of Kent type, by applying a transformation to the von Mises--Fisher density. The results of our test of fit to the von Mises--Fisher law for the subset are displayed in the second row of Table \ref{tab:pal1}. For all significance levels and each choice of the tuning parameter, the tests reject the null hypothesis, indicating a poor fit of the von Mises--Fisher family for the subset of the data.

As a second parametric family of  distributions, we consider the Kent distribution, defined by the density
\begin{equation}\label{eq:denKent}
 f(x)=\frac{1}{c(A,\kappa)} \exp(\kappa \: x^\top \theta+x^\top Ax),\quad x\in\mathcal{S}^{d-1}.
\end{equation}
Here, $\kappa>0$ is a concentration parameter and $\theta\in\mathcal{S}^{d-1}$ is the mean direction.
Moreover, $A$ is a symmetric $d\times d$-matrix with tr$(A)=0$ and $A\theta=0$ that depends on an 'ovality' parameter $\beta$, see \cite{MJ:2000}. Hence we test the hypothesis that the data stems from a density of type \eqref{eq:denKent}, where the parameters $\kappa,\theta$ and $\beta$ are unknown. These parameters have been estimated by the method of maximum-likelihood, and the same bootstrap parameters are applied as above. The bootstrap $p$-values are reported in Table \ref{tab:pal1}. Interestingly, the full data set is rejected by each of the tests on every level of significance. We thus conclude that the Kent distribution is not a suitable model. However, we obtain a different impression  for the subset of the data with age fixed to 1250. For this data set, none of the  tests can reject the hypothesis of an underlying Kent distribution.

\begin{table}[t]
\setlength{\tabcolsep}{.7mm}
    \centering
        \caption{MLE estimates ($\widehat{\kappa}_n$, $\widehat{\theta}_n$) of the parameters of the von Mises--Fisher distribution and bootstrap $p$-values of the test $T_{n,\gamma}^{\{2\}}$ applied to the paleomagnetic data set.}
    \label{tab:pal1}

    \begin{tabular}{l|c|ccc|cccccc}
       model & sample  & $\widehat{\kappa}_n$ & $\widehat{\theta}_n$ & $\widehat{\beta}_n$ & $T_{n,0.1}^{\{2\}}$ & $T_{n,0.5}^{\{2\}}$ & $T_{n,0.75}^{\{2\}}$ & $T_{n,1}^{\{2\}}$ & $T_{n,2}^{\{2\}}$ & $T_{n,3}^{\{2\}}$ \\ \hline
       \multirow{2}{1cm}{vMF} & full sample  & 0.131 & (0.947, -0.319,  0.039) & & 0.99 & 0.024 & 0.264 & 1 & 1 & 1 \\
       & subsample  & 0.587 & (0.282, -0.930,  0.237) & & 0 & 0 & 0 & 0 & 0 & 0 \\\hline
       \multirow{2}{1cm}{Kent} & full sample  & 0.131 & (0.947, -0.319,  0.039) & 0.070 &  0 & 0 & 0 & 0 & 0 & 0 \\
        & subsample  & 0.597 & (0.282, -0.930,  0.237) & 1.043 & 0.206 & 0.376 & 0.287 & 0.19 & 0.291 & 0.419
    \end{tabular}
\end{table}

The results obtained in this section confirm the statements about the data sets made in \cite{SW:2019}.
\section{Discussion} \label{conclude}
 We have studied goodness-of-fit tests for spherical and hyperspherical data. Our tests apply to both simple hypotheses with all parameters assumed known and to composite hypotheses, with parameters estimated from the data at hand. Limit theory is developed under the null hypothesis as well as under alternatives, while the asymptotic validity of a resampling version of the tests is established. The new procedures perform well in finite samples,  and they are competitive against other methods, whenever such methods are available. An application illustrates the usefulness of the new tests for data-modelling on the sphere.

\bibliography{lit_GOF_CFSP}
\bibliographystyle{imsart-number}





\end{document}